\documentclass[a4paper,10pt]{article}
\usepackage{stmaryrd}
\usepackage{amsfonts}
\usepackage{bbm}
\usepackage{amscd}
\usepackage{mathrsfs}
\usepackage{latexsym,amssymb,amsmath,amscd,amscd,amsthm,amsxtra}
\usepackage[dvips]{graphicx}
\usepackage[utf8]{inputenc}
\usepackage[T1]{fontenc}
\usepackage{lmodern}
\usepackage{amssymb}
\usepackage[all]{xy}
\usepackage{nicefrac,mathtools,enumitem}
\usepackage{microtype}
\usepackage{xcolor}
\usepackage{color}
\definecolor{red}{rgb}{1,0,0}
\definecolor{green}{rgb}{0,1,0}
\definecolor{blue}{rgb}{0,0,1}
\definecolor{refkey}{gray}{.625}
\definecolor{labelkey}{gray}{.625}

\newtheorem{thm}{Theorem}[section]
\newtheorem{lem}[thm]{Lemma}
\newtheorem{cor}[thm]{Corollary}
\newtheorem{pro}[thm]{Proposition}

\newtheorem{rmk}[thm]{Remark}
\newtheorem{defi}[thm]{Definition}

\newcommand {\yh}[1]{{\marginpar{*}\scriptsize\textcolor{purple}{yh: #1}}}
\newcommand {\jf}[1]{{\marginpar{*}\scriptsize\textcolor{red}{jf: #1}}}
\newcommand {\emptycomment}[1]{}

\newcommand{\con }{ \mathsf{con}}
\newcommand{\lin }{ \mathsf{lin}}
\newcommand{\cl }{ \mathsf{cl}}

\newcommand{\be }{\begin{equation}}
\newcommand{\ee }{\end{equation}}

\newcommand{\pf}{\noindent{\bf Proof.}\ }

\newcommand{\huaL}{\mathcal{L}}

\newcommand{\huaT}{\mathcal{T}}

\newcommand{\CWM}{C^{\infty}(M)}

\newcommand{\frkl}{\mathfrak l}

\newcommand{\frku}{\mathfrak u}
\newcommand{\frkv}{\mathfrak v}
\newcommand{\frkw}{\mathfrak w}

\newcommand{\frkG}{\mathfrak G}

\newcommand{\frkL}{\mathfrak L}

\newcommand{\frkX}{\mathfrak X}

\def\qed{\hfill ~\vrule height6pt width6pt depth0pt}

\newcommand{\pair}[1]{\left( #1\right)_+}

\newcommand{\Courant}[1]{\left\llbracket  #1\right\rrbracket }

\newcommand{\Dorfman}[1]{\{ #1\}}

\newcommand{\br}[1]{   [ \cdot,    \cdot  ]   }

\newcommand{\g}{\mathfrak g}

\newcommand{\dM}{\mathrm{d}}

\newcommand{\Der}{\mathrm{Der}}

\newcommand{\gl}{\mathfrak {gl}}

\newcommand{\ad}{\mathrm{ad}}

\newcommand{\pr}{\mathrm{pr}}

\begin{document}
\title{Omni $n$-Lie algebras and linearization of higher analogues of Courant algebroids
\thanks
 {
 Research supported by NSFC (11471139), NSF of Jilin Province (20140520054JH,20170101050JC) and Nan Hu Scholar Development Program of XYNU.
 }}

\author{Jiefeng Liu$^1$, Yunhe Sheng$^{1,2}$ and Chunyue Wang$^{2,3}$\\
$^1$Department of Mathematics, Xinyang Normal University,\\ \vspace{2mm} Xinyang 464000, Henan, China\\
 $^2$Department of Mathematics, Jilin University,
 \\\vspace{2mm}Changchun 130012, Jilin, China\\
$^3$ Department of Mathematics, Jilin Engineering Normal University,\\\vspace{2mm} Changchun 130052, Jilin, China
\\ Email:liujf12@126.com; shengyh@jlu.edu.cn; wang1chun2yue3@163.com\\
}

\date{}
\maketitle
\begin{abstract}
In this paper, we introduce the notion of an omni $n$-Lie algebra and show that they are linearization of higher analogues of standard Courant algebroids. We also introduce the notion of a nonabelian omni $n$-Lie algebra and show that they are linearization of higher analogues of   Courant algebroids associated to Nambu-Poisson manifolds.

\end{abstract}

\section{Introduction}

Courant algebroids were introduced in \cite{LWXmani} (see also \cite{Roytenberg4}), and have many applications. See \cite{Schwarzbach4} and references therein for more details.  On $\huaT^{n-1} M\triangleq TM\oplus \wedge^{n-1}T^*M$,   define a symmetric nondegenerate $\wedge^{n-2}T^* M$-valued pairing $(\cdot,\cdot)_+:\huaT^{n-1} M\times \huaT^{n-1} M\longrightarrow \wedge^{n-2}T^* M$  by
\begin{equation}\label{eq:pair}
(X+\alpha,Y+\beta)_+=i_X\beta+i_Y\alpha,\quad\forall X+\alpha,Y+\beta\in\frkX(M)\oplus\Omega^{n-1}(M),
\end{equation}
and define a bracket operation $\Courant{\cdot,\cdot}:\Gamma(\huaT^{n-1} M)\times \Gamma(\huaT^{n-1} M)\longrightarrow \Gamma(\huaT^{n-1} M)$   by
\begin{equation}\label{eq:Dorfman}
\Courant{X+\alpha,Y+\beta}=[X,Y]+L_X\beta-i_Yd\alpha.
\end{equation}
  The quadruple $(TM\oplus \wedge^{n-1}T^*M,(\cdot,\cdot)_+,\Dorfman{\cdot,\cdot},\pr_{TM})$ is called the higher analogue of the standard Courant algebroid. In particular, if $n=2$, we obtain the standard Courant algebroid. Recently, due to applications in multisymplectic geometry, Nambu-Poisson geometry, $L_\infty$-algebra theory and topological field theory, higher analogues of Courant algebroids are widely studied. See \cite{BiSheng1,BouwknegtJ,Grabowski,GS,hagiwara,hull,JV,Zambon} for more details.

The notion of an omni-Lie algebra was introduced by Weinstein in \cite{Alan} to study the linearization of the standard  Courant algebroid. Then it
was studied from several aspects
\cite{kinyon-weinstein,ShengZhu,UchinoOmni}. An {\bf omni-Lie algebra} associated to a vector space $V$ is a triple $(\gl(V)\oplus V,(\cdot,\cdot)_+,\Dorfman{\cdot,\cdot})$, where $(\cdot,\cdot)_+$ is the $V$-valued pairing  given by
\begin{equation}\label{eq:Vpair}
( A+u,B+v)_+=Av+Bu,\quad \forall ~A+u,B+v\in\gl(V)\oplus V,
\end{equation}
and $\Dorfman{\cdot,\cdot}$ is the bilinear bracket operation given by
\begin{equation}\label{omni}
\{A+u,B+v\}=[A,B]+Av.
\end{equation}
 Even though $(\gl(V)\oplus V,\Dorfman{\cdot,\cdot})$ is not a Lie algebra, its Dirac structures characterize all Lie algebra structures on $V$.
 We can construct a Lie 2-algebra from an omni-Lie algebra. See \cite{ShengZhu} for more details.

In \cite{LSXnonabelian}, the authors introduced the notion of a nonabelian omni-Lie algebra $(\gl(\g)\oplus \g,( \cdot,\cdot )_+,\Dorfman{\cdot,\cdot}_\g)$ associated to a Lie algebra $(\g,[\cdot,\cdot]_\g)$, which originally comes from the study of homotopy Poisson manifolds \cite{LangShengXu}. In particular, they showed that it is the linearization of the Courant algebroid $T\g^*\oplus T^*_{\pi_\g}\g^*$ associated to the linear Poisson manifold $(\g^*,\pi_\g)$, where $ \pi_\g$ is the Lie-Poisson structure on $\g^*$.

The purpose of this paper is to extend the above results to the $n$-ary case. First we introduce the notion of an omni $n$-Lie algebra, which is a triple $(\gl(V)\oplus \wedge^{n-1}V,(\cdot,\cdot)_+,\Dorfman{\cdot,\cdot})$ including a bracket operation $\Dorfman{\cdot,\cdot}$ and a $(V\otimes \wedge^{n-2}V)$-valued pairing $(\cdot,\cdot)_+$. Similar to the classical case, $(\gl(V)\oplus \wedge^{n-1}V,\Dorfman{\cdot,\cdot})$ is a Leibniz algebra. We show that a linear map $F:\wedge^{n}V\longrightarrow  V$ defines an $n$-Lie algebra structure on $V$ if and only if the graph of $F^\sharp:\wedge^{n-1}V\longrightarrow\gl(V)$ is a sub-Leibniz algebra of  $(\gl(V)\oplus \wedge^{n-1}V,\Dorfman{\cdot,\cdot})$. Note that this result is not totally parallel the classical case. Namely the condition that $F$ being skew-symmetric can not be simply described by being isotropic with respect to the $(V\otimes \wedge^{n-2}V)$-valued pairing $(\cdot,\cdot)_+$. We further show that an omni $n$-Lie algebra  $(\gl(V)\oplus \wedge^{n-1}V,(\cdot,\cdot)_+,\Dorfman{\cdot,\cdot})$ can be viewed as the linearization of the higher analogue of the standard Courant algebroid $(TM\oplus \wedge^{n-1}T^*M,(\cdot,\cdot)_+,\Courant{\cdot,\cdot},\pr_{TM})$ via letting $M=V^*$. Then we introduce the notion of a nonabelian omni $n$-Lie algebra $(\gl(\g)\oplus \wedge^{n-1}\g,(\cdot,\cdot)_+,\Dorfman{\cdot,\cdot}_\g)$ associated to an $n$-Lie algebra $\g$ and study its algebraic properties. Finally, we give the notion of higher analogues of Courant algebroids associated to Nambu-Poisson manifolds and study their properties. Furthermore, we show that nonabelian omni $n$-Lie algebras are linearization of higher analogues of Courant algebroids associated to Nambu-Poisson manifolds.

The paper is organized as follows. In Section 2, we recall $n$-Lie algebras and Nambu-Poisson manifolds. In Section 3, we introduce the notion of an omni $n$-Lie algebra associated to a vector space $V$ and characterize $n$-Lie algebra structures on $V$ via sub-Leibniz algebra structures of the omni $n$-Lie algebra. In Section 4, we show that an omni $n$-Lie algebra is the linearization of the higher analogue of the standard Courant algebroid. In Section 5, we introduce the notion of a nonabelian omni $n$-Lie algebra and study its algebraic properties. In Section 6, we introduce the notion of higher analogues of Courant algebroids associated to Nambu-Poisson manifolds and show that nonabelian omni $n$-Lie algebras are their linearization.

\section{Preliminaries}

 In this section, we briefly recall the notions of $n$-Lie algebras and Nambu-Poisson manifolds. The notion of  an $n$-Lie algebra, or a Filippov algebra, was introduced in  \cite{Filippov} and have many applications in mathematical physics. See the review article \cite{review} for more details.
 Nambu-Poisson structures were
introduced in \cite{TakhtajanNambu} by Takhtajan in order to find an
axiomatic formalism for Nambu-mechanics which is a generalization of
Hamiltonian mechanics

\begin{defi}\label{defi of n-LA}
An $n$-Lie algebra is a vector space $\g$ together with an $n$-multilinear  skew-symmetric bracket $[\cdot,\cdots,\cdot]_\g:\wedge^n\g\longrightarrow\g$ such that for all $u_i,v_i\in\g$, the following Fundamental  Identity is satisfied:
\begin{eqnarray}\label{FI}
[u_1,u_2,\cdots,u_{n-1},[v_1,v_2,\cdots,v_n]_\g]_\g=\sum_{i=1}^{n}[v_1,v_2,\cdots,[u_1,u_2,\cdots,u_{n-1},v_i]_\g,\cdots,v_n]_\g.
\end{eqnarray}
\end{defi}

For $u_1,u_2,\cdots,u_{n-1}\in\g$, define $\ad:\wedge^{n-1}\g\longrightarrow\gl(\g)$ by $$\ad_{u_1,u_2,\cdots,u_{n-1}}v=[u_1,u_2,\cdots,u_{n-1},v]_\g, \quad\forall ~v\in \g.$$
Then Eq. \eqref{FI} is equivalent to that $\ad_{u_1,u_2,\cdots,u_{n-1}}$ is a derivation, i.e.
  \begin{eqnarray}\label{FI2}
\ad_\frku[v_1,v_2,\cdots,v_n]_\g=\sum_{i=1}^{n}[v_1,v_2,\cdots,\ad_\frku v_i,\cdots,v_n]_\g,\quad \forall~ \frku=u_1\wedge u_2\wedge\cdots\wedge u_{n-1}\in\wedge^{n-1}\g.
\end{eqnarray}
Elements in $\wedge^{n-1}\g$ are called {\bf fundamental objects} of the $n$-Lie algebra $(\g,[\cdot,\cdots,\cdot]_\g)$. In the sequel, we will denote $\ad_\frku v$ simply by $\frku\circ v$.

Define a bilinear  operation on the set of fundamental objects $ \circ: (\wedge^{n-1}\g)\otimes(\wedge^{n-1}\g)\longrightarrow \wedge^{n-1}\g$  by
\begin{eqnarray}\label{comp FOS}
\frku\circ \frkv=\sum_{i=1}^{n-1}v_1\wedge\cdots \wedge v_{i-1}\wedge\frku\circ v_i\wedge v_{i+1}\wedge\cdots\wedge v_{n-1},
\end{eqnarray}
for all $\frku=u_1\wedge u_2\wedge\cdots\wedge u_{n-1}$ and $\frkv=v_1\wedge v_2\wedge \cdots\wedge v_{n-1}.$
In \cite{DT}, the authors proved that $(\wedge^{n-1}\g,\circ)$ is a Leibniz algebra. See \cite{Loday} for details about Leibniz algebras, which are also called Loday algebras. Moreover, the Fundamental Identity (\ref{FI}) is equivalent to
\begin{eqnarray}\label{FI3}
\frku\circ (\frkv\circ w)-\frkv\circ (\frku\circ w)=(\frku\circ \frkv)\circ w,\quad\forall ~\frku,\frkv\in\wedge^{n-1}\g,w\in\g.
\end{eqnarray}

\emptycomment{
A Poisson structure on an $m$-dimensional smooth manifold $M$ is a
bilinear skew-symmetric map
$\{\cdot,\cdot\}:\CWM\times\CWM\longrightarrow\CWM$ such that the
Leibniz rule and the Jacobi identity are satisfied. It is well known
that it is equivalent to a bivector field $\pi\in\wedge^2\frkX(M)$
such that  $[\pi,\pi]=0$. The relation is given by $\{f,g\}=\pi(\dM
f,\dM g)$. We usually denote a Poisson manifold by $(M,\pi)$. For
any $\sigma\in\frkX^{n+1}(M))$ and $\theta\in\Omega^{n+1}(M))$,
$\sigma^\sharp:\wedge^nT^*M\longrightarrow TM$ and
$\theta_\sharp:TM\longrightarrow \wedge^nT^*M$ are given by
$$
\sigma^\sharp(\xi)=i_\xi\sigma,\quad\theta_\sharp(X)=i_X\theta,\quad\forall~\xi\in\Omega^n(M),\forall~X\in\frkX(M).
$$
Denote by $\frkG_\sigma\subset TM\oplus\wedge^nT^*M$ (resp.
$\frkG_\theta$) the graph of $\sigma^\sharp$ (resp.
$\theta_\sharp$).

For a Poisson manifold $(M,\pi)$, there is a skew-symmetric bracket
operation $[\cdot,\cdot]_\pi$ on $1$-forms which is given by
\begin{equation}\label{defi:bracketpi}
[\alpha,\beta]_\pi=L_{\pi^\sharp(\alpha)}\beta-L_{\pi^\sharp(\beta)}\alpha+\dM
i_{\pi^\sharp(\beta)}\alpha,\quad\forall~\alpha,\beta\in\Omega^1(M).
\end{equation}
It follows that $(T^*M,[\cdot,\cdot]_\pi,\pi^\sharp)$ is a Lie
algebroid. Also we have $ [\dM f,\dM g]_\pi=\dM\{f,g\}. $
}

\begin{defi}{\rm\cite{TakhtajanNambu}}
A Nambu-Poisson structure of order $n-1$ on $M$ is an $n$-linear
map
$\{\cdot,\cdots,\cdot\}:\CWM\times\cdots\times\CWM\longrightarrow\CWM$
satisfying the following properties:
\begin{itemize}
\item[\rm(1)] skewsymmetry, i.e. for all $f_1,\cdots,f_{n}\in \CWM$ and $\sigma\in
\mathrm{Sym}(n),$
$$
\{f_1,\cdots,f_{n}\}=(-1)^{|\sigma|}\{f_{\sigma(1)},\cdots,f_{\sigma(n)}\};
$$

\item[\rm(2)]  the Leibniz rule,
i.e. for all $g\in\CWM$, we have
$$
\{f_1g,\cdots,f_{n}\}=f_1\{g,\cdots,f_{n}\}+g\{f_1,\cdots,f_{n}\};
$$

\item[\rm(3)] integrability condition, i.e. for all $f_1,\cdots,f_{n-1},g_1,\cdots,g_n\in \CWM$,
$$\{f_1,\cdots,f_{n-1},\{g_1,\cdots,g_{n}\}\}=\sum_{i=1}^{n}\{g_1,\cdots,\{f_1,\cdots,f_{n-1},g_i\},\cdots,g_{n}\}.$$
\end{itemize}
\end{defi}

In particular, a Nambu-Poisson structure of order $1$ is exactly a
usual Poisson structure. Similar to the fact that a Poisson
structure is determined by a bivector field, a Nambu-Poisson
structure of order $n-1$ is determined by an $n$-vector field
$\pi\in\frkX^{n}(M)$ such that
$$
\{f_1,\cdots,f_{n}\}=\pi(d f_1,\cdots,d f_{n}).
$$
An $n$-vector field $\pi\in\frkX^{n}(M)$ is a Nambu-Poisson
structure if and only if for all $f_1,\cdots,f_{n-1}\in\CWM$, we have
$$
L_{\pi^\sharp(d f_1\wedge\cdots\wedge d f_{n-1})}\pi=0,
$$
where $\pi^\sharp:\wedge^{n-1}T^*M\longrightarrow TM$ is defined by
$$
\langle\pi^\sharp(\xi_1\wedge\cdots\wedge\xi_{n-1}),\xi_n\rangle=\pi(\xi_1\wedge\cdots\wedge\xi_{n-1}\wedge\xi_n),\quad\forall\xi_1,\cdots,\xi_n\in\Omega^1(M).
$$
Let $\pi$ be a Nambu-Poisson structure on a manifold $M$. Then there is a natural operation $[\cdot,\cdot]_\pi$ on $\Omega^{n-1}(M)$ given by
\begin{equation}\label{eq:Nambu-Poisson bracket}
[\alpha,\beta]_\pi=L_{\pi^\sharp(\alpha)}\beta-L_{\pi^\sharp(\beta)}\alpha+d
i_{\pi^\sharp(\beta)}\alpha,\quad\forall~\alpha,\beta\in\Omega^{n-1}(M)
\end{equation}
such that $(\wedge^{n-1}T^*M,[\cdot,\cdot]_\pi,\pi^\sharp)$ is a Leibniz algebroid. See \cite{BiSheng1,Leibnizalgebroid} for more details.

\section{Omni $n$-Lie algebras}

Let $V$ be a finite dimensional vector space. For all $A\in\gl(V)$, define $\huaL_A:\otimes^{n-1}V\longrightarrow \otimes^{n-1}V$ by
$$\huaL_A(v_1\otimes \cdots\otimes v_{n-1})=\sum_{i=1}^{n-1}v_1\otimes\cdots\otimes Av_i \otimes \cdots\otimes v_{n-1}.$$

\begin{defi}
An {\bf omni $n$-Lie algebra} associated to a vector space $V$ is a triple $(\gl(V)\oplus \wedge^{n-1} V,(\cdot,\cdot)_+,\{\cdot,\cdot\})$, where $\{\cdot,\cdot\}$ is the bilinear bracket operation given by
\begin{equation}\label{eq:omnin2}
\{A+\frku,B+\frkv\}=[A,B]+\huaL_A\frkv,
\end{equation}
and $(\cdot,\cdot)_+$ is the $(V\otimes \wedge^{n-2}V)$-valued pairing given by
\begin{equation}\label{eq:omnin1}
(A+\frku,B+\frkv)_+=\sum^{n-1}_{i=1}(-1)^{i+1}\big(Av_i\otimes v_1\wedge\cdots \wedge\hat{v}_i\wedge\cdots\wedge v_{n-1}+Bu_i\otimes u_1\wedge\cdots \wedge\hat{u}_i\wedge\cdots \wedge u_{n-1}\big),
\end{equation}
where $\frku=u_1\wedge u_2\wedge\cdots\wedge u_{n-1}$ and $\frkv=v_1\wedge v_2\wedge\cdots\wedge v_{n-1}$.
\end{defi}

\begin{rmk}
  When $n=2$, we recover Weinstein's omni-Lie algebras \cite{Alan}.
\end{rmk}

\begin{pro}
With the above notations, $(\gl(V)\oplus \wedge^{n-1} V,\{\cdot,\cdot\})$ is a Leibniz algebra. Furthermore, the pairing $(\cdot,\cdot)_+$ and the bracket $\{\cdot,\cdot\}$ are compatible in the sense that
\begin{equation}\label{eq:omnin3}
(\{e_1,e_2\},e_3)_++(e_2,\{e_1,e_3\})_+=\rho_V(e_1)(e_2,e_3)_+,
\end{equation}
where $e_i\in\gl (V)\oplus \wedge^{n-1}V,~i=1,2,3$, and $\rho_V:\gl(V)\oplus \wedge^{n-1} V\longrightarrow \gl(V\otimes\wedge^{n-2} V)$ is given by
\begin{equation}
\rho_V(A+\frku)(\frkw)=\huaL_A\frkw,\quad \forall~\frkw\in V\otimes\wedge^{n-2} V.
\end{equation}
\end{pro}
\pf Since $\huaL_{[A,B]}=\huaL_A\huaL_B-\huaL_B\huaL_A$, we can deduce that $(\gl(V)\oplus \wedge^{n-1} V,\{\cdot,\cdot\})$ is a Leibniz algebra directly.

For all $A+\frku,B+\frkv,C+\frkw\in\gl(V)\oplus \wedge^{n-1} V$, on one hand, we have
\begin{eqnarray*}
&&( \{A+\frku,B+\frkv\},C+\frkw )_++( B+\frkv,\{A+\frku,C+\frkw\})_+\\ &=&([A,B]+\huaL_A\frkv,C+\frkw)_++(B+\frkv,[A,C]+\huaL_A\frkw)_+\\
&=&\sum_{j=1}^{n-1}(-1)^{j+1}\big(ABw_j\otimes w_1\wedge\cdots \wedge\hat{w}_j\wedge\cdots\wedge w_{n-1}+ACv_j\otimes v_1\wedge\cdots \wedge\hat{v}_j\wedge\cdots\wedge v_{n-1}\big)\\
&&+\sum_{i\neq j}^{n-1}(-1)^{i+1}\big(Bw_i\otimes w_1\wedge\cdots \wedge{Aw_j}\wedge\cdots \wedge w_{n-1})+Cv_i\otimes v_1\wedge\cdots \wedge{Av_j}\wedge\cdots \wedge v_{n-1}\big).
\end{eqnarray*}
On the other hand, we have
\begin{eqnarray*}
&&\rho_V(A+u)( B+v, C+w)_+\\
&=&\sum^{n-1}_{j=1}(-1)^{j+1}\rho_V(A+u)\big(Bw_j\otimes w_1\wedge\cdots \wedge\hat{w}_j\wedge\cdots\wedge w_{n-1}+Cv_j\otimes v_1\wedge\cdots \wedge \hat{v}_j\wedge\cdots\wedge v_{n-1}\big)\\
&=&\sum_{j=1}^{n-1}(-1)^{j+1}\big(ABw_j\otimes w_1\wedge\cdots \wedge\hat{w}_j\wedge\cdots\wedge w_{n-1}+ACv_j\otimes v_1\wedge\cdots \wedge\hat{v}_j\wedge\cdots\wedge v_{n-1}\big)\\
&&+\sum_{i\neq j}^{n-1}(-1)^{i+1}\big(Bw_i\otimes w_1\wedge\cdots \wedge{Aw_j}\wedge\cdots \wedge w_{n-1}+Cv_i\otimes v_1\wedge\cdots \wedge{Av_j}\wedge\cdots \wedge v_{n-1}\big),
\end{eqnarray*}
which implies that \eqref{eq:omnin3} holds. \qed\vspace{3mm}

Let $F:\wedge^nV\longrightarrow V$ be a linear map. Then $F$ induces a linear map $F^\sharp:\wedge^{n-1}V\longrightarrow \gl(V)$ by
$$
F^\sharp(\frku)(u)=F(\frku,u),\quad\forall \frku\in \wedge^{n-1}V,~u\in V.
$$
Denote by $G_F\subset \gl(V)\oplus\wedge^{n-1}V$ the graph of $F^\sharp$.

\begin{thm}
Let $F:\wedge^nV\longrightarrow V$ be a linear map. Then $(V,F)$ is an $n$-Lie algebra if and only if $G_F$ is a Leibniz subalgebra of the Leibniz algebra $(\gl(V)\oplus \wedge^{n-1} V,\{\cdot,\cdot\})$.
\end{thm}
\pf $G_F$ is a Leibniz subalgebra of the Leibniz algebra $(\gl(V)\oplus \wedge^{n-1} V,\{\cdot,\cdot\})$ if and only if for all $\frku,\frkv\in\wedge^{n-1} V$, $\{F^\sharp(\frku)+\frku,F^\sharp(\frkv)+\frkv\}\in G_F$, which is equivalent to
$$F^\sharp(\huaL_{F^\sharp(\frku)}\frkv)=[F^\sharp(\frku),F^\sharp(\frkv)].$$
Since $\huaL_{F^\sharp(\frku)}\frkv=\sum_{i=1}^{n-1}v_1\wedge\cdots\wedge F(\frku,v_i)\wedge\cdots\wedge v_{n-1}$, thus the above equality can be written as
$$F^\sharp(\frku\circ\frkv)=[F^\sharp(\frku),F^\sharp(\frkv)],$$
which is equivalent to that $(V,F)$ is an $n$-Lie algebra.\qed

\section{Linearization of the higher analogue of the standard Courant algebroid}\label{sec:LHSCA}
Let $V$ be an $m$-dimensional vector space and $V^*$ its dual space. We consider the direct sum bundle $\huaT^{n-1}V^*=TV^*\oplus\wedge^{n-1}T^*V^*$. Denote the vector spaces of linear vector fields and constant $(n-1)$-forms on $V^*$ by $\frkX_{\lin}(V^*)$ and $\Omega^{n-1}_{\con}(V^*)$ respectively.  It is obvious that $\frkX_{\lin}(V^*)\oplus\Omega^{n-1}_{\con}(V^*)\cong\gl(V)\oplus\wedge^{n-1}V$. To make this explicit, for any $x\in V$, denote by $l_x$ the corresponding linear function on $V^*$. Let $\{x^i\}$ be a basis of the vector space $V$. Then $\{l_{x^i}\}$ forms a coordinate chart for $V^*$. So $\{\frac{\partial}{\partial l_{x^i}}\}$ constitutes a basis of vector fields on $V^*$ and $\{d l_{x^i}\}$ constitutes a basis of $1$-forms on $V^*$.
For $A=(a_{ij})\in\gl(V)$, we get a linear vector field $\widehat{A}=\sum_il_{A(x^i)}\frac{\partial}{\partial l_{x^i}}$ on $V^*$. Also $\frku=\sum_{1\leq i_1< \cdots <i_{n-1}\leq m}\xi_{i_1\cdots i_{n-1}}x^{i_1}\wedge\cdots\wedge x^{i_{n-1}}$ defines a constant ${(n-1)}$-form $\widehat{\frku}=\sum_{1\leq i_1< \cdots <i_{n-1}\leq m}\xi_{i_1\cdots i_{n-1}}dl_{x^{i_1}}\wedge\cdots\wedge dl_{x^{i_{n-1}}}$.

  Define $\Phi:\gl(V)\oplus\wedge^{n-1}V\longrightarrow \frkX_{\lin}(V^*)\oplus\Omega^{n-1}_{\con}(V^*)$ by
  $$\Phi(A+\frku)=\widehat{A}+\widehat{\frku}.$$
Obviously, $\Phi$ is an isomorphism between vector spaces.

Any element $v\otimes\frku\in V\otimes \wedge^{n-2}V$ will give rise to a linear $(n-2)$-form $\overline{v\otimes\frku}$ defined by
$$ \overline{v\otimes\frku}=l_{v}\widehat{\frku}.$$

We give some useful formulas first.
\begin{lem}\label{lem:LHSCA}
With the above notations, for all $A\in\gl(V)$ and $\frku\in\wedge^{n-1}V$, we have
\begin{eqnarray}
(\widehat{A},\widehat{\frku})_+&=&\overline{(A,\frku)}_+\label{eq:form1},\\
d(\widehat{A},\widehat{\frku})_+&=&\widehat{\huaL_A\frku}\label{eq:form2},\\
L_{\widehat{A}}\widehat{\frku}&=&\widehat{\huaL_A\frku}\label{eq:form3},\\
{[\widehat{A},\widehat{B}]}&=&\widehat{[A,B]}\label{eq:form4}.
\end{eqnarray}
\end{lem}
\pf On one hand, for $\frku=\sum_{1\leq i_1< \cdots <i_{n-1}\leq m}\xi_{i_1\cdots i_n}x^{i_1}\wedge\cdots\wedge x^{i_{n-1}}\in \wedge^{n-1}V$, we have
\begin{eqnarray*}
(\widehat{A},\widehat{\frku})_+=\sum_{1\leq i_1< \cdots <i_{n-1}\leq m}\sum_{j=1}^n(-1)^{j+1} \xi_{i_1\cdots i_j\cdots i_n}l_{A{x^{i_j}}}dl_{x^{i_1}}\wedge\cdots\wedge \hat{dl_{x^{i_j}}}\wedge\cdots\wedge dl_{x^{i_{n-1}}}.
\end{eqnarray*}
On the other hand, we have
\begin{eqnarray*}
\overline{(A,\frku)}_+&=&\sum_{1\leq i_1< \cdots <i_{n-1}\leq m}\sum_{j=1}^{n-1}(-1)^{j+1}\xi_{i_1\cdots i_n}\overline{(A{x^{i_j}}\otimes x^{i_1}\wedge\cdots\wedge\hat{x^{i_j}}\wedge \cdots\wedge x^{i_{n-1}})}\\
&=&\sum_{1\leq i_1< \cdots <i_{n-1}\leq m}\sum_{j=1}^{n-1}(-1)^{j+1} \xi_{i_1\cdots i_j\cdots i_{n-1}}l_{A{x^{i_j}}}dl_{x^{i_1}}\wedge\cdots\wedge \hat{dl_{x^{i_j}}}\wedge\cdots\wedge dl_{x^{i_{n-1}}},
\end{eqnarray*}
which implies that \eqref{eq:form1} holds.

By direct calculation, we have
\begin{eqnarray*}
d(\widehat{A},\widehat{\frku})_+&=&\sum_{1\leq i_1< \cdots <i_n\leq m}\sum_{j=1}^{n-1}(-1)^{j+1} \xi_{i_1\cdots i_j\cdots i_n}dl_{Ax^{i_j}}\wedge dl_{x^{i_1}}\wedge\cdots\wedge \hat{dl_{x^{i_j}}}\wedge\cdots\wedge dl_{x^{i_{n-1}}}\\
&=&\sum_{1\leq i_1< \cdots <i_{n-1}\leq m}\sum_{j=1}^n\sum_{k=1}^m(-1)^{j-1}a_{k,i_j}\xi_{i_1\cdots i_j\cdots i_{n-1}}d l_{x^k}\wedge dl_{x^{i_1}}\wedge\cdots\wedge \hat{dl_{x^{j}}}\wedge\cdots\wedge dl_{x^{i_{n-1}}}.\\
\end{eqnarray*}
On the other hand, we have
\begin{eqnarray*}
\huaL_A\frku&=&A(\sum_{1\leq i_1< \cdots <i_{n-1}\leq m}\xi_{i_1\cdots i_p}{x^{i_1}}\wedge\cdots\wedge {x^{i_{n-1}}})\\
&=&\sum_{1\leq i_1< \cdots <i_{n-1}\leq m}\sum_{j=1}^{n-1}\xi_{i_1\cdots i_{n-1}}{x^{i_1}}\wedge\cdots\wedge Ax^{i_j}\wedge\cdots\wedge {x^{i_{n-1}}}\\
&=&\sum_{1\leq i_1< \cdots <i_{n-1}\leq m}\sum_{j=1}^{n-1}\sum_{k=1}^ma_{k,i_j}\xi_{i_1\cdots i_{n-1}}{x^{i_1}}\wedge\cdots\wedge x^k\wedge\cdots\wedge {x^{i_{n-1}}}.
\end{eqnarray*}
Thus \eqref{eq:form2} follows immediately.

\eqref{eq:form3} follows from
$$L_{\widehat{A}}\widehat{u}=\iota_{\widehat{A}}d\widehat{u}+d\iota_{\widehat{A}}\widehat{u}=d(\widehat{A},\widehat{\frku})_+=\widehat{\huaL_Au}.$$

\eqref{eq:form4} is direct. We omit the details. \qed\vspace{3mm}

Now we are ready to show that an omni $n$-Lie algebra can be viewed as linearization of the higher analogue of the standard Courant algebroid.

\begin{thm}\label{thm:LHSCA}
The omni $n$-Lie algebra $(\gl(V)\oplus \wedge^{n-1} V,(\cdot,\cdot)_+,\{\cdot,\cdot\})$ is induced from the higher analogue of the standard Courant algebroid $(\huaT^{n-1}(V^*),\pair{\cdot,\cdot},\Courant{\cdot,\cdot},\pr_{TV^*})$ via restriction to $\frkX_{\lin}(V^*)\oplus\Omega^{n-1}_{\con}(V^*)$. More precisely, we have
\begin{eqnarray}
\label{eq:ind1} ( \Phi({A}+ {\frku}),\Phi({B}+{\frkv}))_+&=&\overline{(A+\frku,B+\frkv )}_+,\\
\label{eq:ind2} \Courant{\Phi (A+ {\frku}),\Phi({B}+{\frkv})}&=&\Phi{\{A+\frku,B+\frkv\}},\\
\label{eq:ind3} L_{\pr_{TV^*}\Phi(A+\frku)}\overline{\frkw}&=&\overline{\rho_V(A+\frku)(\frkw)}.
\end{eqnarray}
\end{thm}
\pf By \eqref{eq:form1}, we have
\begin{eqnarray*}
 (\Phi({A}+{\frku}),\Phi({B}+{\frkv}))_+&=&(\widehat{A},\widehat{\frkv})_++(\widehat{\frku},\widehat{B})_+=\overline{(A,\frku)}_++\overline{(\frkv,B)}_+=\overline{(A+\frku,B+\frkv)}_+.
\end{eqnarray*}
By \eqref{eq:form3} and \eqref{eq:form4}, we have
\begin{eqnarray*}
\Courant{\Phi (A+ {\frku}),\Phi({B}+{\frkv})}&=&\Courant{\widehat{A}+\widehat{\frku},\widehat{B}+\widehat{\frkv}}=[\widehat{A},\widehat{B}]+L_{\widehat{A}}\widehat{\frkv}\\
&=&\widehat{[A,B]}+\widehat{\huaL_A\frkv}=\Phi{\{A+\frku,B+\frkv\}}.
\end{eqnarray*}
Finally, for all $\frkw=w_1\otimes w_2\wedge \cdots\wedge w_{n-1}\in V\otimes \wedge^{n-2}V$, we have
\begin{eqnarray*}
L_{\pr_{TV^*}\Phi(A+\frku)}\overline{\frkw}&=&L_{\widehat{A}}(l_{w_1}dw_2\wedge\cdots\wedge dw_{n-1})\\
&=&L_{\widehat{A}}(l_{w_1})dw_2\wedge\cdots\wedge dw_{n-1}+l_{w_1}(\sum_{i=2}^{n-1}dw_2\wedge\cdots \wedge L_{\widehat{A}}dw_i\wedge\cdots\wedge dw_{n-1})\\
&=&l_{Aw_1}dw_2\wedge\cdots\wedge dw_{n-1}+\sum_{i=2}^{n-1}l_{w_1}dw_2\wedge\cdots \wedge d(Aw_i)\wedge\cdots\wedge dw_{n-1})\\
&=&\overline{Aw_1\otimes w_2\wedge \cdots\wedge w_{n-1}}+\sum_{i=2}^{n-1}\overline{{w_1}\otimes w_2\wedge\cdots \wedge Aw_i\wedge\cdots\wedge w_{n-1}}\\
&=&\overline{\rho_V(A+\frku)(\frkw)}.
\end{eqnarray*}
The proof is finished.\qed

\emptycomment{\begin{cor}
The relation between the compatible conditions of the higher analogues of the standard Courant algebroid $(\huaT^{n-1}(V^*),\pair{\cdot,\cdot},\Courant{\cdot,\cdot},\rho)$ and  the omni $n$-Lie algebra $(\gl(V)\oplus \wedge^{n-1} V,(\cdot,\cdot)_+,\{\cdot,\cdot\})$ is given by
\begin{eqnarray*}
(\Courant{\hat{e}_1,\hat{e}_2},\hat{e}_3)_++(\hat{e}_2,\Courant{\hat{e}_1,\hat{e}_3})_+-\rho(\hat{e}_1)(\hat{e}_2,\hat{e}_3)_+=\frkl_{((\{e_1,e_2\},e_3)_++(e_2,\{e_1,e_3\})_+-\rho_V(e_1)(e_2,e_3)_+)}.
\end{eqnarray*}
\end{cor}
\yh{Both hands are zero. What do you mean?}}

\section{Nonabelian omni $n$-Lie algebras}

\begin{defi}
A {\bf nonabelian omni $n$-Lie algebra} associated to an $n$-Lie algebra $(\g,[\cdot,\cdots,\cdot]_\g)$ is a triple $(\gl(\g)\oplus \wedge^{n-1} \g,(\cdot,\cdot)_+,\{\cdot,\cdot\}_\g)$, where $(\cdot,\cdot)_+$ is the $(\g\otimes \wedge^{n-1}\g)$-valued pairing given by \eqref{eq:omnin1} and $\{\cdot,\cdot\}_\g$ is the bilinear bracket operation given by
\begin{equation}\label{eq:omninn2}
\{A+\frku,B+\frkv\}_\g=[A,B]+[A,\ad_\frkv]+[\ad_\frku,B]-\ad_{\huaL_A\frkv}+\huaL_A\frkv+\frku\circ\frkv,\quad\forall ~A,B\in\gl(\g),~\frku,\frkv\in\wedge^{n-1}\g.
\end{equation}
\emptycomment{The following bracket is also a Leibniz algebra.
\begin{equation}\label{eq:omnin22}
\{A+\frku,B+\frkv\}_\g=[A,B]+[\ad_\frku,B]+\huaL_A\frkv+\frku\circ\frkv,
\end{equation}
The condition $[A,\ad_\frkv]=\ad_{\huaL_A\frkv}$ means that $A$ is a derivation for $n$-Lie algebra $(\g,[\cdot,\cdots,\cdot]_\g).$  which one is more meaningful? It seems that \eqref{eq:omnin2} is more like the one for nonabelian omni-Lie algebras. For \eqref{eq:omnin22}, you need that $A$ is a derivation for $n$-Lie algebra $(\g,[\cdot,\cdots,\cdot]_\g)$? Why, in general $A$ do not need to be a derivation. Thus, I guess you are considering a subalgebra $\Der(\g)\oplus \wedge^{n-1}\g$}

\end{defi}
\begin{thm}
 Let $(\gl(\g)\oplus \wedge^{n-1} \g,(\cdot,\cdot)_+,\{\cdot,\cdot\}_\g)$ be a nonabelian omni $n$-Lie algebra. Then we have
 \begin{itemize}
   \item[\rm(i)] $(\gl(\g)\oplus \wedge^{n-1} \g,\{\cdot,\cdot\}_\g)$ is a Leibniz algebra;

    \item[\rm(ii)] $\{A+\frku,A+\frku\}_\g=-\ad_{\huaL_A\frku}+\huaL_A\frku+\frku\circ\frku$;

     \item[\rm(iii)]the pairing $(\cdot,\cdot)_+$ and the bracket $\{\cdot,\cdot\}_\g$   are compatible in the sense that
\begin{equation}\label{eq:omninn3}
\rho_\g(e_1)(e_2,e_3)_+=(\{e_1,e_2\}_\g-([A,\ad_\frkv]-\ad_{\huaL_A\frkv}),e_3)_++(e_2,\{e_1,e_3\}_\g-([A,\ad_\frkw]-\ad_{\huaL_A\frkw}))_+,
\end{equation}
where $e_1=A+\frku, e_2=B+\frkv, e_3=C+\frkw\in\gl (\g)\oplus \wedge^{n-1}\g$ and $\rho_\g:\gl(\g)\oplus \wedge^{n-1} \g\longrightarrow \gl(\g\otimes\wedge^{n-2} \g)$ is given by
\begin{equation}
\rho_\g(A+\frku)(\frkw)=\huaL_{A+\ad_{\frku}}\frkw,\quad \forall ~\frkw\in \g\otimes\wedge^{n-2} \g.
\end{equation}
 \end{itemize}
\end{thm}
\pf (i) We can prove that $(\gl(\g)\oplus \wedge^{n-1} \g,\{\cdot,\cdot\}_\g)$ is a Leibniz algebra directly by a complicated computation. In the sequel, we will show that $(\gl(\g)\oplus \wedge^{n-1} \g,\{\cdot,\cdot\}_\g)$ is a trivial deformation of the omni $n$-Lie algebra $(\gl(\g)\oplus \wedge^{n-1} \g,\{\cdot,\cdot\})$. Thus, we omit details here.

(ii) It follows from \eqref{eq:omninn2} directly.

(iii) By straightforward computation, we have
\begin{eqnarray*}
&&([A,B]+[\ad_{\frku},B]+\huaL_A\frkv+\frku\circ\frkv,C+\frkw)_+\\
&=&\sum^{n-1}_{i=1}(-1)^{i+1}([A,B]+[\ad_{\frku},B])w_i\otimes w_1\wedge\cdots \wedge\hat{w_i}\wedge\cdots\wedge w_{n-1}\\
&&+\sum_{i \neq j}(-1)^{j+1}Cv_{j}\otimes v_1\wedge \cdots \wedge\hat{v_j}\wedge\cdots \wedge (Av_i+\frku\circ v_i)\wedge\cdots\wedge v_{n-1}\\
&&+\sum^{n-1}_{i=1}(-1)^{i+1}(CAv_i+C(\frku\circ v_i))\otimes v_1\wedge \cdots  \wedge \widehat{v_i}\wedge\cdots\wedge v_{n-1}.
\end{eqnarray*}
On the other hand, we have
\begin{eqnarray*}
&&(B+\frkv,[A,C]+[\ad_{\frku},C]+\huaL_A\frkw+\frku\circ\frkw)_+\\
&=&\sum^{n-1}_{i=1}(-1)^{i+1}([A,C]+[\ad_{\frku},C])v_i\otimes v_1\wedge\cdots \wedge\widehat{v_i}\wedge\cdots\wedge v_{n-1}\\
&&+\sum_{i \neq j}(-1)^{j+1}(Bw_{j}\otimes w_1\wedge \cdots \wedge\hat{w_j}\wedge\cdots \wedge (Aw_i+\frku\circ w_i)\wedge\cdots\wedge w_{n-1})\\
&&+\sum^{n-1}_{i=1}(-1)^{i+1}(BAw_i+B(\frku\circ w_j))\otimes w_1\wedge \cdots  \wedge \hat{ w_i}\wedge\cdots\wedge w_{n-1}.
\end{eqnarray*}
Thus we have
\begin{eqnarray*}
&&([A,B]+[\ad_{\frku},B]+\huaL_A\frkv+\frku\circ\frkv,C+\frkw)_++(B+\frkv,[A,C]+[\ad_{\frku},C]+\huaL_A\frkw+\frku\circ\frkw)_+\\
&=&\sum^{n-1}_{i=1}(-1)^{i+1}(A\circ B+\ad_{\frku}\circ B)w_i\otimes w_1\wedge\cdots \wedge\hat{w_i}\wedge\cdots\wedge w_{n-1}\\
&&+\sum^{n-1}_{i=1}(-1)^{i+1}(A\circ C+\ad_{\frku}\circ C)v_i\otimes v_1\wedge\cdots \wedge\hat{v_i}\wedge\cdots\wedge v_{n-1}\\
&&+\sum_{i \neq j}(-1)^{j+1}Bv_{j}\otimes v_1\wedge \cdots \wedge\hat{v_j}\wedge\cdots \wedge (Av_i+{\frku}\circ v_i)\wedge\cdots\wedge v_{n-1}\\
&&+\sum_{i \neq j}(-1)^{j+1}Bw_{j}\otimes w_1\wedge \cdots \wedge\hat{w_j}\wedge\cdots \wedge (Aw_i+{\frku}\circ w_i)\wedge\cdots\wedge w_{n-1}\\
&=&\rho_\g(A+\frku)(B+\frkv,C+\frkw)_+.
\end{eqnarray*}
The proof is finished. \qed\vspace{3mm}

Obviously, for all $A\in\Der(\g)$, we have $[A,\ad_\frkv]-\ad_{\huaL_A\frkv}=0$.
Thus, we have

\begin{cor}
  For all $e_1,e_2,e_3\in \Der(\g)\oplus \wedge^{n-1}\g$, we have
  $$
  \rho_\g(e_1)(e_2,e_3)_+=(\{e_1,e_2\}_\g,e_3)_++(e_2,\{e_1,e_3\}_\g)_+.
  $$
\end{cor}

 In the sequel we show that a nonabelian omni $n$-Lie algebra $(\gl(\g)\oplus \wedge^{n-1} \g,\{\cdot,\cdot\}_\g)$ can be viewed as a trivial deformation of the omni $n$-Lie algebra $(\gl(\g)\oplus \wedge^{n-1} \g,\{\cdot,\cdot\})$. For details of deformations of Leibniz algebras, see \cite{Car,Kos}.

Let $(\frkL,[\cdot,\cdot]_\frkL)$ be a Leibniz algebra. For an endomorphism $N$ of $\frkL$, define
$$[e_1,e_2]_N=[Ne_1,e_2]_\frkL+[e_1,Ne_2]_\frkL-N[e_1,e_2]_\frkL,$$
and set
$$TN(e_1,e_2)=[Ne_1,Ne_2]_\frkL-N[e_1,e_2]_N.$$
The endomorphism $N$ is called a {\bf Nijenhuis operator}  if $TN=0$.

A Nijenhuis operator gives  a trivial  deformation of the Leibniz algebra $(\frkL,[\cdot,\cdot]_\frkL)$.
\begin{pro}{\rm\cite{Car}}\label{pro:N}
Let $N$ be a Nijenhuis operator on the Leibniz algebra $(\frkL,[\cdot,\cdot]_\frkL)$. Then we have
\begin{itemize}
  \item[\rm(1)]  $(\frkL,[\cdot,\cdot]_N)$ is a Leibniz algebra;
    \item[\rm(2)] $N$ is a morphism of Leibniz algebras from $(\frkL,[\cdot,\cdot]_N)$ to $(\frkL,[\cdot,\cdot]_\frkL)$;
      \item[\rm(3)] $(\frkL,[\cdot,\cdot]_\frkL+[\cdot,\cdot]_N)$  is a Leibniz algebra.
\end{itemize}
\end{pro}
Let  $(\g,[\cdot,\cdots,\cdot]_\g)$ be an $n$-Lie algebra. Then we can define a linear map $N:\gl(\g)\oplus \wedge^{n-1}\g\longrightarrow \gl(\g)\oplus \wedge^{n-1}\g$ by
\begin{equation}\label{eq:N}
N(A+\frku)=\ad_\frku.
\end{equation}

\begin{lem}\label{lem:N}
The linear map $N$ given by \eqref{eq:N} is a Nijenhuis operator on the Leibniz algebra $(\gl(\g)\oplus \wedge^{n-1}\g,\{\cdot,\cdot\})$, where the Leibniz bracket $\{\cdot,\cdot\}$ is given by $(\ref{eq:omnin2})$.
\end{lem}
\pf First by definition, we have
\begin{eqnarray*}
\{A+\frku,B+\frkv\}_N&=&\{N(A+\frku),B+\frkv\}+\{A+\frku,N(B+\frkv)\}-N\{A+\frku,B+\frkv\}\\ &=&[\ad_\frku,B]+\huaL_{\ad_\frku}\frkv+[A,\ad_\frkv]-\ad_{\huaL_A\frkv}\\
&=&[\ad_\frku,B]+\frku\circ\frkv+[A,\ad_\frkv]-\ad_{\huaL_A\frkv}.
\end{eqnarray*}
Hence it is clear that
$$N\{A+\frku,B+\frkv\}_N=\ad_{\frku\circ\frkv}=[\ad_\frku,\ad_\frkv]=\{N(A+\frku),N(B+\frkv)\},$$
which says that $N$ is a Nijenhuis operator.\qed\vspace{3mm}

It is straightforward  to see that
$$\{A+\frku,B+\frkv\}_\g=\{A+\frku,B+\frkv\}+\{A+\frku,B+\frkv\}_N.$$
Therefore, by Proposition \ref{pro:N} and Lemma \ref{lem:N}, we have
\begin{thm}
 Let  $(\g,[\cdot,\cdots,\cdot]_\g)$ be an $n$-Lie algebra. Then the bracket $\{\cdot,\cdot\}_\g$ is a trivial deformation of the Leibniz bracket $\{\cdot,\cdot\}$. In particular, $(\gl(\g)\oplus\wedge^{n-1}\g, \{\cdot,\cdot\}_\g)$ is a Leibniz algebra.
\end{thm}

\begin{rmk}
 If we view $(\gl(\g),[\cdot,\cdot])$ as a Leibniz algebra, then $(\gl(\g),[\cdot,\cdot])$ and $(\wedge^{n-1}\g,\circ)$ form a matched pair of Leibniz algebras and the Leibniz algebra $(\gl(\g)\oplus\wedge^{n-1}\g,\{\cdot,\cdot\}_\g)$ is exactly their double. See \cite{Agore} for more details about matched pairs of Leibniz algebras.
\end{rmk}

\section{Linearization of higher analogues of Courant algebroids associated to Nambu-Poisson structures}
Let $(M,\pi)$ be a Nambu-Poisson   manifold. We introduce a bracket $\Courant{\cdot,\cdot}_\pi:\Gamma(\huaT^{n-1}M)\times\Gamma(\huaT^{n-1}M)\longrightarrow \Gamma(\huaT^{n-1}M)$  by
\begin{equation}\label{eq:Nambu-Courant Bracket}
\Courant{X+\alpha,Y+\beta}_\pi=[X,Y]+[X,\pi^\sharp(\beta)]+[\pi^\sharp(\alpha),Y]-\pi^\sharp(L_{X}\beta)+\pi^\sharp(i_Yd\alpha)+L_{X}\beta-i_Yd\alpha+[\alpha,\beta]_\pi,
\end{equation}
where $X,Y\in\frkX(M),\alpha,\beta\in\Omega^{n-1}(M)$ and $[\cdot,\cdot]_\pi$ is given by \eqref{eq:Nambu-Poisson bracket}.

Let $\rho_\pi:\huaT^{n-1}M\rightarrow TM$ be the bundle map defined by
\begin{equation}\label{eq:Nambu-Courant anchor}
\rho_\pi(X+\alpha)=X+\pi^\sharp(\alpha),\quad\forall\ X\in\frkX(M),\alpha\in\Omega^{n-1}(M).
\end{equation}
We call the quadruple $(\huaT^{n-1}M,(\cdot,\cdot)_+,\Courant{\cdot,\cdot}_\pi,\rho_\pi)$ {\bf the higher analogue of the Courant algebroid associated to a Nambu-Poisson manifold} and denote it by $\huaT^{n-1}_\pi M$. In the sequel, we will see that even though we call it the higher analogue of a Courant algebroid, some important properties for Courant algebroids do not hold anymore.

\begin{thm}
Let $(\huaT^{n-1}M,(\cdot,\cdot)_+,\Courant{\cdot,\cdot}_\pi,\rho_\pi)$ be the higher analogue of the Courant algebroid associated to a Nambu-Poisson manifold. Then we have
\begin{itemize}
\item[\rm(i)]   $(\huaT^{n-1}M,\Courant{\cdot,\cdot}_\pi,\rho_\pi)$ is a
Leibniz algebroid.
\item[\rm(ii)] The bracket $\Courant{\cdot,\cdot}_\pi$ is not skew-symmetric. Instead, we have
\begin{equation}\label{eq:skew}
\Courant{X+\alpha,X+\alpha}_\pi= d(X,\alpha)_++d(\pi^\sharp(\alpha),\alpha)_+-\pi^\sharp(d(X,\alpha)_+).
\end{equation}
\item[\rm(iii)]The pairing \eqref{eq:pair} and the bracket $\Courant{\cdot,\cdot}_\pi$ are
compatible in the following sense:
\begin{eqnarray*}\label{eq:relation}
L_{\rho(e_1)}(e_2,e_3)_+&=&\pair{\Courant{e_1,e_2}_\pi-([X,\pi^\sharp(\beta)]-\pi^\sharp(L_X\beta)+\pi^\sharp(i_Yd\alpha))+i_{\pi^\sharp(\beta)}d\alpha,e_3}\\
&&+\pair{e_2,\Courant{e_1,e_3}_\pi-([X,\pi^\sharp(\gamma)]-\pi^\sharp(L_X\gamma)+\pi^\sharp(i_Zd\alpha))+i_{\pi^\sharp(\gamma)}d\alpha}.
\end{eqnarray*}
where $e_1=X+\alpha,~e_2=Y+\beta,~e_3=Z+\gamma$.
\end{itemize}
\end{thm}

\pf (i) Let $\Psi:\huaT^{n-1}M\longrightarrow \huaT^{n-1}M$ be the invertible bundle map defined by
\begin{equation}\label{eq:isom}
\Psi(X+\alpha)=X+\alpha+\pi^\sharp(\alpha),\quad\forall\ X\in\frkX(M),\alpha\in\Omega^{n-1}(M).
\end{equation}
By direct calculation, we have
\begin{equation}
\Psi\Courant{X+\alpha,Y+\beta}_\pi=\Courant{\Psi(X+\alpha),\Psi(Y+\beta)},
\end{equation}
where the bracket $\Courant{\cdot,\cdot}$ is given by \eqref{eq:Dorfman}. Since $(\Gamma(\huaT^{n-1}M),\Courant{\cdot,\cdot})$ is a Leibniz algebra, we deduce that $(\Gamma(\huaT^{n-1}M),\Courant{\cdot,\cdot}_\pi)$ is also a Leibniz algebra.

 For all $f\in \CWM$,   $X,Y\in\frkX(M)$
and $\alpha,\beta\in\Omega^{n-1}(M)$,   we have
\begin{eqnarray*}
\Courant{X+\alpha,f(Y+\beta)}_\pi&=&[X,fY]+[X,f\pi^\sharp(\beta)]+[\pi^\sharp(\alpha),fY]-\pi^\sharp(L_{X}f\beta)\\
&&+\pi^\sharp(i_{fY}d\alpha)+L_{X}f\beta-i_{fY}d\alpha+[\alpha,f\beta]_\pi\\
&=&f[X,Y]+X(f)Y+f[X,\pi^\sharp(\beta)]+X(f)\pi^\sharp(\beta)+f[\pi^\sharp(\alpha),Y]\\
&&+\pi^\sharp(\alpha)(f)Y-f\pi^\sharp(L_{X}\beta)-X(f)\pi^\sharp(\beta)+f\pi^\sharp(i_{Y}d\alpha)+fL_X\beta\\
&&+X(f)\beta-fi_{Y}d\alpha+f[\alpha,\beta]_\pi+\pi^\sharp(\alpha)(f)\beta\\
&=&f\Courant{X+\alpha,Y+\beta}_\pi+(X+\pi^\sharp(\alpha))(f)(Y+\beta).
\end{eqnarray*}

Thus, $ ( \huaT^{n-1}M),\Courant{\cdot,\cdot}_\pi,\rho_\pi)$ is a Leibniz algebroid.

(ii) It is straightforward to obtain \eqref{eq:skew} by \eqref{eq:Nambu-Courant Bracket}.

(iii) The left hand side of the above equality is equal to
$$L_Xi_Y\gamma+L_Xi_Z\beta+L_{\pi^\sharp(\alpha)}i_Y\gamma+L_{\pi^\sharp(\alpha)}i_Z\beta.$$
The right hand side is equal to
\begin{eqnarray*}
&&i_ZL_X\beta+i_ZL_{\pi^\sharp(\alpha)}\beta-i_ZL_{\pi^\sharp(\beta)}\alpha+i_Zdi_{\pi^\sharp(\beta)}\alpha+i_{[X,Y]}\gamma\\
&&+i_{[X,\pi^\sharp(\beta)]}\gamma+i_{[\pi^\sharp(\alpha),Y]}\gamma-i_{\pi^\sharp(L_X\beta)}\gamma+i_{\pi^\sharp(i_Yd\alpha)}\gamma\\
&&+i_YL_X\gamma+i_YL_{\pi^\sharp(\alpha)}\gamma-i_YL_{\pi^\sharp(\gamma)}\alpha+i_Ydi_{\pi^\sharp(\gamma)}\alpha+i_{[X,Z]}\beta\\
&&+i_{[X,\pi^\sharp(\gamma)]}\beta +i_{[\pi^\sharp(\alpha),Z]}\beta-i_{\pi^\sharp(L_X\gamma)}\beta+i_{\pi^\sharp(i_Zd\alpha)}\beta.
\end{eqnarray*}
The conclusion follows from
$$i_{[X,Y]}=L_Xi_Y-i_YL_X.$$
This finishes the proof.\qed\vspace{3mm}

Let $\frkX_H(M)$ and $\Omega_{\cl}^{n-1}(M)$ denote the set of the Hamiltonian vector fields and closed $(n-1)$-forms respectively.

\begin{cor}
For all $e_1, ~e_2,~e_3\in\frkX_H(M)\oplus \Omega_{\cl}^{n-1}(M)$, we have
\begin{equation}\label{eq:Nambu-Courant compatible2}
L_{\rho_\pi(e_1)}(e_2,e_3)_+=\pair{\Courant{e_1,e_2}_\pi,e_3}+\pair{e_2,\Courant{e_1,e_3}_\pi}.
\end{equation}
\end{cor}
\pf For all $e_1=X+\alpha,~e_2=Y+\beta,~e_3=Z+\gamma\in\frkX_H(M)\oplus \Omega_{\cl}^{n-1}(M)$, since $\alpha$ is closed, we have
$$i_{\pi^\sharp(i_Yd\alpha)}\gamma-i_Yi_{\pi^\sharp(\gamma)}d\alpha+i_{\pi^\sharp(i_Zd\alpha)}\beta-i_Zi_{\pi^\sharp(\beta)}d\alpha=0.$$
For all $\xi,\eta\in\Omega^{n-1}(M)$, we have the following formula
$$\pi^\sharp(L_{\pi^\sharp(\xi)}\eta)-[\pi^\sharp(\xi),\pi^\sharp(\eta)]=(-1)^{n-1}(i_{d\xi}\pi)\pi^\sharp(\eta).$$
Without loss of generality, let $X=\pi^\sharp(d f_1\wedge\cdots \wedge df_{n-1})$, then we have
\begin{eqnarray*}
&&i_{[X,\pi^\sharp(\beta)]}\gamma-i_{\pi^\sharp(L_X\beta)}\gamma+i_{[X,\pi^\sharp(\gamma)]}\beta-i_{\pi^\sharp(L_X\gamma)}\beta\\
&=&(-1)^ni_{(i_{d(d f_1\wedge\cdots \wedge df_{n-1})}\pi)\pi^\sharp(\beta)}\gamma+(-1)^ni_{(i_{d(d f_1\wedge\cdots \wedge df_{n-1})}\pi)\pi^\sharp(\gamma)}\beta=0.
\end{eqnarray*}
We finishes the proof.\qed
\vspace{3mm}

In the following, we show that the nonabelian omni $n$-Lie algebra is a linearization of the higher analogue of the Courant algebroid $(\huaT^{n-1}M,(\cdot,\cdot)_+,\Courant{\cdot,\cdot}_\pi,\rho_\pi)$ associated to a Nambu-Poisson manifold $(M,\pi)$.

Let $(\g,[\cdot,\cdots,\cdot]_\g)$ be an $m$-dimensional $n$-Lie algebra such that it induces a linear Nambu-Poisson structure\footnote{Not all $n$-Lie algebras can give linear Nambu-Poisson structures on dual spaces, see \cite{DufZung} for details.} $\pi_\g$ on $\g^*$.
 Then we obtain the higher analogue of the Courant algebroid
  $\huaT^{n-1}_{\pi_\g}\g^*$. Let $\{x^i\}$ be a basis of the vector space $\g$. Using the same notations as in Section \ref{sec:LHSCA},  we have
$$\pi_\g=\sum_{1\leq i_1< \cdots <i_{n}\leq m}l_{[x^{i_1},\cdots,x^{i_n}]_{\g}}\frac{\partial}{\partial l_{x^{i_1}}}\wedge \cdots\wedge\frac{\partial}{\partial l_{x^{i_n}}}.$$

  \emptycomment{Denote the vector spaces of linear vector fields and constant $(n-1)$-forms on $\g^*$ by $\frkX_{\lin}(\g^*)$ and $\Omega^{n-1}_{\con}(\g^*)$ respectively.  It is obvious that $\frkX_{\lin}(\g^*)\oplus\Omega^{n-1}_{\con}(\g^*)\cong\gl(\g)\oplus\wedge^{n-1}\g$. To make this explicit, for any $x\in \g$, denote by $l_x$ the corresponding linear function on $\g^*$. Let $\{x^i\}$ be a basis of the vector space $\g$. Then $\{l_{x^i}\}$ forms a coordinate chart for $\g^*$. So $\{\frac{\partial}{\partial l_{x^i}}\}$ constitutes a basis of vector fields on $\g^*$ and $\{d l_{x^i}\}$ constitutes a basis of $1$-forms on $\g^*$.
For $A\in\gl(\g)$, we get a linear vector field $\widehat{A}=\sum_jl_{A(x^j)}\frac{\partial}{\partial l_{x^j}}$ on $\g^*$. Also $\frku=\sum_{1\leq i_1< \cdots <i_{n-1}\leq m}\xi_{i_1\cdots i_{n-1}}x^{i_1}\wedge\cdots\wedge x^{i_{n-1}}$ defines a constant ${n-1}$-form $\widehat{\frku}=\sum_{1\leq i_1< \cdots <i_{n-1}\leq m}\xi_{i_1\cdots i_{n-1}}dl_{x^{i_1}}\wedge\cdots\wedge dl_{x^{i_{n-1}}}$. Moreover,
}

\begin{lem}
For all $A\in\gl(\g)$ and $\frku,\frkv\in\wedge^{n-1}\g$, we have
\begin{eqnarray}
\label{eq:t2}\pi_\g^\sharp (\widehat{\frku})&=&\widehat{\ad_\frku},\\
\label{eq:t4}{[\widehat{\frku},\widehat{\frkv}]}_{\pi_\g}&=&\widehat{\frku\circ\frkv},\\
\label{eq:t6}\pi^\sharp(L_{\widehat{A}}\widehat{\frku})&=&\widehat{\ad_{\huaL_A\frku}}.
\end{eqnarray}
\end{lem}
\pf For $\frku=\sum_{1\leq j_1< \cdots <j_{n-1}\leq m}\xi_{j_1\cdots j_{n-1}}x^{j_1}\wedge\cdots\wedge x^{j_{n-1}}\in\wedge^{n-1}\g$ with the corresponding constant ${(n-1)}$-form $\widehat{\frku}=\sum_{1\leq j_1< \cdots <j_{n-1}\leq m}\xi_{j_1\cdots j_{n-1}}dl_{x^{j_1}}\wedge\cdots\wedge dl_{x^{j_{n-1}}}$, we have
\begin{eqnarray*}
\pi^\sharp_\g(\widehat{\frku})&=&\sum_{1\leq i_1< \cdots <i_{n}\leq m}\sum_{k=1}^n(-1)^{n-k}\xi_{i_1\cdots\hat{i_k}\cdots i_n}l_{[x^{i_1},\cdots,x^{i_k},\cdots,x^{i_n}]_\g}\frac{\partial}{\partial l_{x^{i_k}}}\\
&=&\sum_{1\leq i_1< \cdots <i_{n}\leq m}\sum_{k=1}^n\xi_{i_1\cdots\hat{i_k}\cdots i_n}l_{[x^{i_1},\cdots,x^{i_n},x^{i_k}]_\g}\frac{\partial}{\partial l_{x^{i_k}}}\\
&=&\sum_{1\leq j_1< \cdots <j_{n-1}\leq m}\sum_{l=1}^m\xi_{j_1\cdots j_{n-1}}l_{[x^{j_1},\cdots,x^{j_{n-1}},x^{l}]_\g}\frac{\partial}{\partial l_{x^{l}}}\\
&=&\widehat{\ad_\frku},
\end{eqnarray*}
which implies that \eqref{eq:t2} holds.

Since $\widehat{\frku}$ is a constant $(n-1)$-form, by \eqref{eq:form3} and \eqref{eq:t2}, we have
\begin{eqnarray*}
{[\widehat{\frku},\widehat{\frkv}]}_{\pi_\g}&=&L_{\pi^\sharp(\widehat{\frku})}\widehat{\frkv}-L_{\pi^\sharp(\widehat{\frkv})}\widehat{\frku}+d
i_{\pi^\sharp(\widehat{\frkv})}\widehat{\frku}\\
&=&L_{\pi^\sharp(\widehat{\frku})}\widehat{\frkv}-i_{\pi^\sharp(\widehat{\frkv})}d\widehat{\frku}\\
&=&L_{\pi^\sharp(\widehat{\frku})}\widehat{\frkv}=L_{\widehat{\ad_\frku}}\widehat{\frkv}\\
&=&\widehat{\huaL_{\ad_\frku}\frkv}=\widehat{\frku\circ\frkv},
\end{eqnarray*}
which implies that \eqref{eq:t4} holds.

By \eqref{eq:form3} and \eqref{eq:t2}, we have
$$\pi^\sharp(L_{\widehat{A}}\widehat{\frku})=\pi^\sharp(\widehat{\huaL_{A}\frku})=\widehat{\ad_{\huaL_A\frku}}.$$
This ends the proof.\qed

\begin{thm}Let $(\g,[\cdot,\cdots,\cdot]_\g)$ be an $m$-dimensional $n$-Lie algebra such that it induces a linear Nambu-Poisson structure  $\pi_\g$ on $\g^*$.
Then the nonabelian omni $n$-Lie algebra $(\gl(\g)\oplus \wedge^{n-1} \g,(\cdot,\cdot)_+,\{\cdot,\cdot\}_\g)$  is induced from the higher analogue of the  Courant algebroid $(\huaT^{n-1}\g^*,(\cdot,\cdot)_+,\Courant{\cdot,\cdot}_{\pi_\g},\rho_{\pi_\g})$   associated to the Nambu-Poisson manifold $(\g^*,\pi_\g)$ via restriction to $\frkX_{\lin}(\g^*)\oplus\Omega^{n-1}_{\con}(\g^*)$. More precisely, we have
\begin{eqnarray}
\label{eq:ind12} ( \Phi({A}+ {\frku}),\Phi({B}+{\frkv}))_+&=&\overline{(A+\frku,B+\frkv )}_+,\\
\label{eq:ind22} \Courant{\Phi (A+ {\frku}),\Phi({B}+{\frkv})}_{\pi_\g}&=&\Phi{\{A+\frku,B+\frkv\}}_\g,\\
\label{eq:ind32} L_{\rho_{\pi_\g}\Phi(A+\frku)}\overline{\frkw}&=&\overline{\rho_\g(A+\frku)(\frkw)}.
\end{eqnarray}
\end{thm}
\pf \eqref{eq:ind12} has been proved in Theorem \ref{thm:LHSCA}.
By \eqref{eq:form2}-\eqref{eq:form4} and \eqref{eq:t2}-\eqref{eq:t6}, we deduce that
\begin{eqnarray*}
\Courant{\widehat{A}+ \widehat{\frku},\widehat{B}+\widehat{\frkv}}_{\pi_\g}&=&[\widehat{A},\widehat{B}]+[\widehat{A},\pi^\sharp(\widehat{\frkv})]+[\pi^\sharp(\widehat{\frku}),\widehat{B}]-\pi^\sharp(L_{\widehat{A}}\widehat{\frkv})+\pi^\sharp(i_{\widehat{B}}d\widehat{\frku})+L_{\widehat{A}}\widehat{\frkv}-i_{\widehat{B}}d\widehat{\frku}+[\widehat{\frku},\widehat{\frkv}]_\pi\\
&=&[\widehat{A},\widehat{B}]+[\widehat{A},\widehat{\ad_\frkv}]+[\widehat{\ad_\frku},\widehat{B}]-\widehat{\ad_{\huaL_A\frkv}}+\widehat{\huaL_A\frkv}+\widehat{\frku\circ\frkv}\\
&=&\widehat{[{A},{B}]}+\widehat{[{A},\ad_\frkv]}+\widehat{[\ad_\frku,B]}-\widehat{\ad_{\huaL_A\frkv}}+\widehat{\huaL_A\frkv}+\widehat{\frku\circ\frkv}\\
&=&\Phi{\{A+\frku,B+\frkv\}}_\g,
\end{eqnarray*}
which implies that \eqref{eq:ind22} holds.
By \eqref{eq:form4}, we have $L_{\widehat{A}}\overline{\frkw}=\overline{\rho_\g(A)(\frkw)}$. Thus
\begin{eqnarray*}
L_{\rho_\pi\Phi(A+\frku)}\overline{\frkw}&=&L_{\widehat{A}+\pi^\sharp(\widehat{\frku})}\overline{\frkw}=L_{\widehat{A}}\overline{\frkw}+L_{\widehat{\ad_\frku}}\overline{\frkw}\\
&=&\overline{\rho_\g(A)(\frkw)}+\overline{\rho_\g(\ad_\frku)(\frkw)}=\overline{\rho_\g(A+\frku)(\frkw)},
\end{eqnarray*}
which says that \eqref{eq:ind32} holds. This ends the proof.\qed

{}

\begin{thebibliography}{}
\bibitem{Agore}
A. L. Agore and G. Militaru, Unified products for Leibniz algebras. Applications., \emph{Linear Algebra  Appl.} 2013, 439(9): 2609-2633.



\bibitem{BiSheng1}
Y. Bi and Y. Sheng, On higher analogues of Courant algebroids, \emph{Sci. China Math.}  Vol. 54 (3) (2011),   437-447.



\bibitem{BouwknegtJ}
P. Bouwknegt and B. Jurcč, AKSZ construction of topological open p-brane action and Nambu brackets, \emph{Rev. Math. Phys.} 25 (2013), 1330004, 31 pages.


\bibitem{Car}
 J. Cariñena, J. Grabowski and G. Marmo, Courant algebroid and Lie bialgebroid
contractions, \emph{J. Phys. A: Math. Gen.} 2004, 37(19): 5189-5202.


 \bibitem{DT} Y. Daletskii and L. Takhtajan, Leibniz and Lie algebra structures for Nambu
algebra, \emph{Lett. Math. Phys.} 39 (1997), 127-141.

 \bibitem{review}
J. A. de Azc$\rm\acute{a}$rraga and J. M. Izquierdo, $n$-ary algebras: a review with applications,
\emph{ J. Phys. A: Math. Theor.} 43 (2010), 293001.

\bibitem{DufZung}
J. P. Dufour and N. T. Zung, Linearization of Nambu structures,
\emph{Compositio Math.} 117 (1999), no. 1, 77–98.


\bibitem{Filippov}  V. T. Filippov, $n$-Lie algebras,  {\it Sib. Mat. Zh.} \textbf{26} (1985) 126-140.


\bibitem{Grabowski}
J. Grabowski, Brackets,  \emph{Int. J. Geom. Methods Mod. Phys.} 10 (2013), 1360001, 45 pages.

\bibitem{GS}
M. Grützmann and T. Strobl,   General Yang-Mills type gauge theories for p-form gauge fields: from physics-based ideas to a mathematical framework or from Bianchi identities to twisted Courant algebroids, \emph{Int. J. Geom. Methods Mod. Phys.} 12 (2015), no. 1, 1550009, 80 pp.

\bibitem{hagiwara}
Y. Hagiwara, Nambu–Dirac manifolds, \emph{J. Phys. A} 35(5) (2002) 1263-1281.

\bibitem{hull}
C. M. Hull, Generalised geometry for M-theory, \emph{J. High Energy Phys.} 07 (2007) 079.


\bibitem{Leibnizalgebroid}
R. Ibáñez, M. de León, J. Marrero and E. Padrón, Leibniz algebroid associated with a Nambu-
Poisson structure, \emph{J. Phys. A} 32 (46): 8129-8144, 1999

\bibitem{JV}
B. Jurčo and J. Vysoký, Leibniz algebroids, generalized Bismut connections and Einstein–Hilbert actions, \emph{J. Geom. Phys.} 97 (2015),  25-33.

\bibitem{kinyon-weinstein}
M.~K. Kinyon and A.~Weinstein,
\newblock Leibniz algebras, {C}ourant algebroids, and multiplications on
  reductive homogeneous spaces,
\newblock {\em Amer. J. Math.} 123 (3): 525-550, 2001.

\bibitem{Kos}
Y. Kosmann-Schwarzbach, Nijenhuis structures on Courant algebroids, \emph{Bull. Braz. Math. Soc. (N.S.)} 2011, 42(4): 625-649.

\bibitem{Schwarzbach4}
Y. Kosmann-Schwarzbach, Courant algebroids. A short history, \emph{SIGMA Symmetry Integrability Geom. Methods Appl.} 9 (2013), Paper 014, 8 pp.


\bibitem{LangShengXu}
H. Lang, Y. Sheng and X. Xu, Strong homotopy Lie algebras, homotopy Poisson manifolds and Courant algebroids,  \emph{Lett. Math. Phys.} (2016), doi:10.1007/s11005-016-0925-8.

\bibitem{LSXnonabelian}H. Lang, Y. Sheng and X. Xu, Nonabelian omni-Lie algebras, \emph{Banach Center Publications}, 110 (2016), 167-176.

\bibitem{LWXmani}%
Z. Liu, A. Weinstein and P. Xu, Manin triples for Lie
bialgebroids, \emph{J. Diff. Geom.} 1997, 45(3): 547-574.

\bibitem{Loday} J. L. Loday and T. Pirashvilo, Universal enveloping algebras of Leibniz algebras and (co)homology, {\it Math. Ann.} 1993, 296 (1): 139-158.



\bibitem{Roytenberg4}
D. Roytenberg, \emph{Courant algebroids, derived brackets and even
symplectic supermanifolds}, PhD thesis, UC Berkeley, 1999,
arXiv: math.DG/9910078.




\bibitem{ShengZhu}
Y. Sheng and C. Zhu, Semidirect products of representations up to homotopy, {\it Pacific J. Math.} 2011, 249(1): 211-236.

\bibitem{TakhtajanNambu}
L. Takhtajan, On foundation of the generalized {N}ambu mechanics, \emph{ Comm. Math. Phys.} 160(2):295-315, 1994

\bibitem{UchinoOmni}
K.~Uchino.
\newblock Courant brackets on noncommutative algebras and omni-{L}ie algebras,
\newblock {\em Tokyo J. Math.} 30(1):239-255, 2007.

\bibitem{Alan}
A. Weinstein, Omni-Lie algebras, Microlocal analysis of the Schrödinger equation and related topics (Japanese) (Kyoto, 1999), \emph{S${\bar{u}}$rikaisekikenky${\bar{u}}$sho K${\bar{u}}$ky${\bar{u}}$roku}, 2000, 1176: 95-102.


\bibitem{Zambon}
M. Zambon, $L_\infty$-algebras and higher analogues of Dirac structures and Courant
algebroids, \emph{J. Symplectic Geom.} 10(4) (2012), 563-599.
\end{thebibliography}
\end{document}